\documentclass{IEEEtran}

\pdfoutput=1

\usepackage{amsmath}
\usepackage{array}
\usepackage{mdwmath}
\usepackage{eqparbox}                                                                                                                           
\usepackage{float}
\usepackage[font=footnotesize]{subfig}
\usepackage{fixltx2e}
\usepackage{stfloats}
\usepackage{amsfonts}
\usepackage{cite}
\usepackage{amsmath}
\usepackage{bbm}
\usepackage{amssymb}
\usepackage{mathrsfs}
\usepackage{amsbsy}	
\usepackage{color}
\usepackage{multicol}

\usepackage{graphicx,epstopdf}

\usepackage{algorithm}
\usepackage{algpseudocode}
\renewcommand{\algorithmicrequire}{{\bf Input:}}
\renewcommand{\algorithmicensure}{{\bf Output:}}

\usepackage{verbatim}

\floatstyle{ruled}
\newfloat{algorithm}{tbp}{loa}
\floatname{algorithm}{Algorithm}

\long\def\hide#1{}

\usepackage{comment}

\newtheorem{theorem}{\bf Theorem}[section]

\newtheorem{remark}{\bf Remark}[section]

\newcommand{\droped}[1]{{\color{blue} \sout{}}}

\def\ba{\begin{array}}
\def\ea{\end{array}}
\newcommand{\beq}{\begin{equation}}
\newcommand{\eeq}{\end{equation}}
\newcommand{\bq}{\begin{eqnarray}}
\newcommand{\eq}{\end{eqnarray}}
\newcommand{\bqn}{\begin{eqnarray*}}
\newcommand{\eqn}{\end{eqnarray*}}
\newcommand{\bee}{\begin{enumerate}}
\newcommand{\eee}{\end{enumerate}}
\newcommand{\bi}{\begin{itemize}}
\newcommand{\ei}{\end{itemize}}

\newcommand{\mathN}{\mathcal{N}}
\newcommand{\mathE}{\mathcal{E}}
\newcommand{\mathG}{\mathcal{G}}
\newcommand{\mathT}{\mathcal{T}}
\newcommand{\mathI}{\mathcal{I}}
\newcommand{\mathK}{\mathcal{K}}

\newcommand{\mathP}{\mathcal{P}}

\newboolean{showcomments}
\setboolean{showcomments}{true}
\newcommand{\qiuyu}[1]{  \ifthenelse{\boolean{showcomments}}
{ \textcolor{blue}{(Qiuyu says:  #1)}} {}  }
\newcommand{\slow}[1]{\ifthenelse{\boolean{showcomments}}
{ \textcolor{red}{(Steven says:  #1)}}{}}
\newcommand{\pmark}[1]{\ifthenelse{\boolean{showcomments}}
{ \textcolor{cyan}{(Label:  #1)}}{}}

\begin{document}

\title{Distributed Algorithm for Optimal Power Flow on Unbalanced Multiphase Distribution Networks}

\author{Qiuyu Peng and Steven H. Low
\thanks{*This work was supported by ARPA-E grant DE-AR0000226, Los Alamos National Lab through an DoE grant DE-AC52-06NA25396,
DTRA through grant HDTRA 1-15-1-0003 and Skotech.}
\thanks{Qiuyu Peng is with the Electrical Engineering Department 
and Steven H. Low is with the Computing and Mathematical Sciences and
the Electrical Engineering Departments, California Institute of Technology, Pasadena, CA 91125, USA.
{\small \tt \{qpeng, slow\}@caltech.edu}}%
}

\maketitle
\begin{abstract}

The optimal power flow (OPF) problem is fundamental in power distribution networks control and
operation that underlies many important applications such as volt/var control and demand response, etc..
Large-scale highly volatile renewable penetration in the distribution networks
calls for real-time feedback control, and hence the need for distributed solutions for the OPF problem.
Distribution networks are inherently unbalanced and most of the existing distributed solutions for 
balanced networks do not apply. 
In this paper we propose a solution for unbalanced distribution networks.
Our distributed algorithm is based on alternating direction method of 
multiplier (ADMM). Unlike existing approaches that require to solve semidefinite programming problems in each ADMM macro-iteration, we exploit the problem structures and decompose the OPF problem in such a way that the subproblems in each ADMM macro-iteration reduce to either a closed form solution or eigen-decomposition of a $6\times 6$ hermitian matrix, which significantly reduce the convergence time. We present simulations on IEEE 13, 34, 37 and 123 bus unbalanced distribution networks to illustrate the scalability and optimality of the proposed algorithm.
\end{abstract}


\section{Introduction}

The optimal power flow (OPF) problem seeks to minimize a certain objective, such as power loss and 
generation cost subject to power flow physical laws and operational constraints. 
It is a fundamental problem that underpins many distribution system operations and planning problems such as economic dispatch, unit commitment, state estimation, volt/var control and demand response, etc.. Most algorithms proposed in the literature are centralized and meant for applications in today's energy management systems that centrally schedule a relatively small number of generators.   
The increasing penetrations of highly volatile renewable energy sources in distribution systems requires simultaneously optimizing (possibly in real-time) the operation of a large number of intelligent endpoints. A centralized approach will not scale because of its computation and communication overhead and we need to rely on distributed solutions.

Various distributed algorithms for OPF problem have been proposed in the literature. Some early distributed algorithms, including \cite{kim1997coarse,baldick1999fast}, do not deal with the non-convexity issue of OPF and convergence is not guaranteed for those algorithms. Recently, convex relaxation has been applied to convexify the OPF problem, e.g. semi-definite programming (SDP) relaxation \cite{Bai2008,lavaei2012zero,zhang2011geometry,gan2014convex} and second order cone programming (SOCP) relaxation \cite{Jabr2006,Farivar-2013-BFM-TPS,gan2015exact}. When an optimal solution of the original OPF problem can be recovered from any optimal solution of the SOCP/SDP relaxation, we say the relaxation is {exact}. It is shown that both SOCP and SDP relaxations are exact for radial networks using standard IEEE test networks and many practical networks \cite{lavaei2012zero,dall2013distributed,Farivar-2013-BFM-TPS,gan2014convex}. This is important because almost all distribution
systems are radial. Thus, optimization decompositions can be applied to the relaxed OPF problem with guaranteed convergence, e.g. dual decomposition method \cite{lam2012distributed,lam2012optimal}, methods of multiplier \cite{devane2013stability,li2012demand}, and alternating direction method of multiplier (ADMM) \cite{dall2013distributed,kraning2013dynamic,sun2013fully}.

There are at least two challenges in designing distributed algorithm that solves the OPF problem on distribution systems. First, distribution systems are inherently unbalanced because of the unequal loads on each phase \cite{kersting2012distribution}. Most of the existing approaches \cite{lam2012distributed,lam2012optimal,devane2013stability,li2012demand,kraning2013dynamic,sun2013fully,peng2014distributed} are designed for balanced networks and do not apply to unbalanced networks.  

Second, the convexified OPF problem on unbalanced networks consists of semi-definite constraints. To our best knowledge, all the existing distributed solutions \cite{dall2013distributed,kim1997coarse,baldick1999fast} require solving SDPs within each macro-iteration. The SDPs are computationally intensive to solve, and those existing algorithms take significant long time to converge even for moderate size networks. 

In this paper, we address those two challenges through developing an \emph{efficient} distributed algorithm for the OPF problems on \emph{unbalanced} networks based on alternating direction methods of multiplier(ADMM). The advantages of the proposed algorithm are twofold: 1) instead of relying on SDP optimization solver to solve the optimization subproblems in each iteration as existing approaches, we exploit the problem structures and decompose the problem in such a way that the subproblems in each ADMM macro-iteration reduce to either a closed form or a eigen-decomposition of a $6\times 6$ hermitian matrix, which greatly speed up the convergence time. 2) Communication is only required between adjacent buses. 

We demonstrate the scalability of the proposed algorithms using standard IEEE test networks \cite{kersting1991radial}. The proposed algorithm converges within $3$ seconds on the IEEE-13, 34, 37, 123 bus systems. To show the superiority of using the proposed procedure to solve each subproblem, we also compare the computation time for solving a subproblem by our algorithm and an off-the-shelf SDP optimization solver (CVX, \cite{grant2008cvx}). 
Our solver requires on average $3.8\times 10^{-3}$s while CVX requires on average $0.58$s.

A preliminary version has appeared in \cite{peng2015distributed}. In this paper, we improve the algorithm in \cite{peng2015distributed} in the following aspects: 1) We consider more general forms of objective function and power injection region such that the algorithm can be used in more applications. In particular, we provide a sufficient condition, which holds in practice, for the existence of efficient solutions to the optimization subproblems. 2) Voltage magnitude constraints, which are crucial to distribution system operations, are considered in the new algorithm. 3) We study the impact of network topologies on the rate of convergence in the simulations.

The rest of the paper is structured as follows. The OPF problem on an unbalanced network is defined in section \ref{sec:model}. 
In section \ref{sec:alg}, we develop our distributed algorithm based on ADMM.  In section \ref{sec:case}, we
 test its scalability on standard IEEE distribution systems and study the impact of network topologies on the rate of convergence.
  We conclude this paper in section \ref{sec:conclusion}.

\section{Problem Formulation}\label{sec:model}

In this section, we define the optimal power flow (OPF) problem on unbalanced radial distribution networks and review how to solve it through SDP relaxation.

We denote the set of complex numbers with $\mathbb{C}$, the set of $n$-dimensional complex numbers with $\mathbb{C}^n$ and the set of $m\times n$ complex matrix with $\mathbb{C}^{m\times n}$. The set of hermitian (positive semidefinite) matrix is denoted by $\mathbb{S}$ ($\mathbb{S}_+$). The hermitian transpose of a vector (matrix) $x$ is denoted by $x^H$.

The trace of a square matrix $x\in \mathbb{C}^{n\times n}$ is denoted by $tr(x):=\sum_{i=1}^n x_{ii}$. The inner product of two matrices (vectors) $x,y\in\mathbb{C}^{m\times n}$ is denoted by $\langle x,y\rangle:=\mathbf{Re}(tr(x^Hy))$. The Frobenius (Euclidean) norm of a matrix (vector) $x\in\mathbb{C}^{m\times n}$ is defined as $\|x\|_2:=\sqrt{\langle x,x\rangle}$. Given $x\in\mathbb{C}^{n\times n}$, let diag$(x)\in \mathbb{C}^{n\times 1}$ denote the vector composed of $x$'s diagonal elements. 

\subsection{Branch flow model}

We model a distribution network by a \emph{directed} tree graph $\mathT := (\mathN, \mathE)$ where
 $\mathN:=\{0,\ldots,n\}$ represents the set of buses and $\mathE$ represents the set
 of distribution lines connecting the buses in $\mathN$.  
Index the root of the tree by $0$ and let $\mathN_+:=\mathN\setminus\{0\}$ denote the other buses.  For each bus $i$, it has a unique ancestor $A_i$ and a set of children buses, denoted by $C_i$. We adopt the graph orientation where every line points towards the root. Each directed line connects a bus $i$ and its unique ancestor $A_i$. We hence label the lines by $\mathE:=\{1,\ldots,n\}$ where each $i\in\mathE$ denotes a line from $i$ to $A_i$. Note that $\mathE=\mathN_+$ and we will use $\mathN_+$ to represent the lines set for convenience.

\begin{figure}
\centering
\includegraphics[scale=0.25]{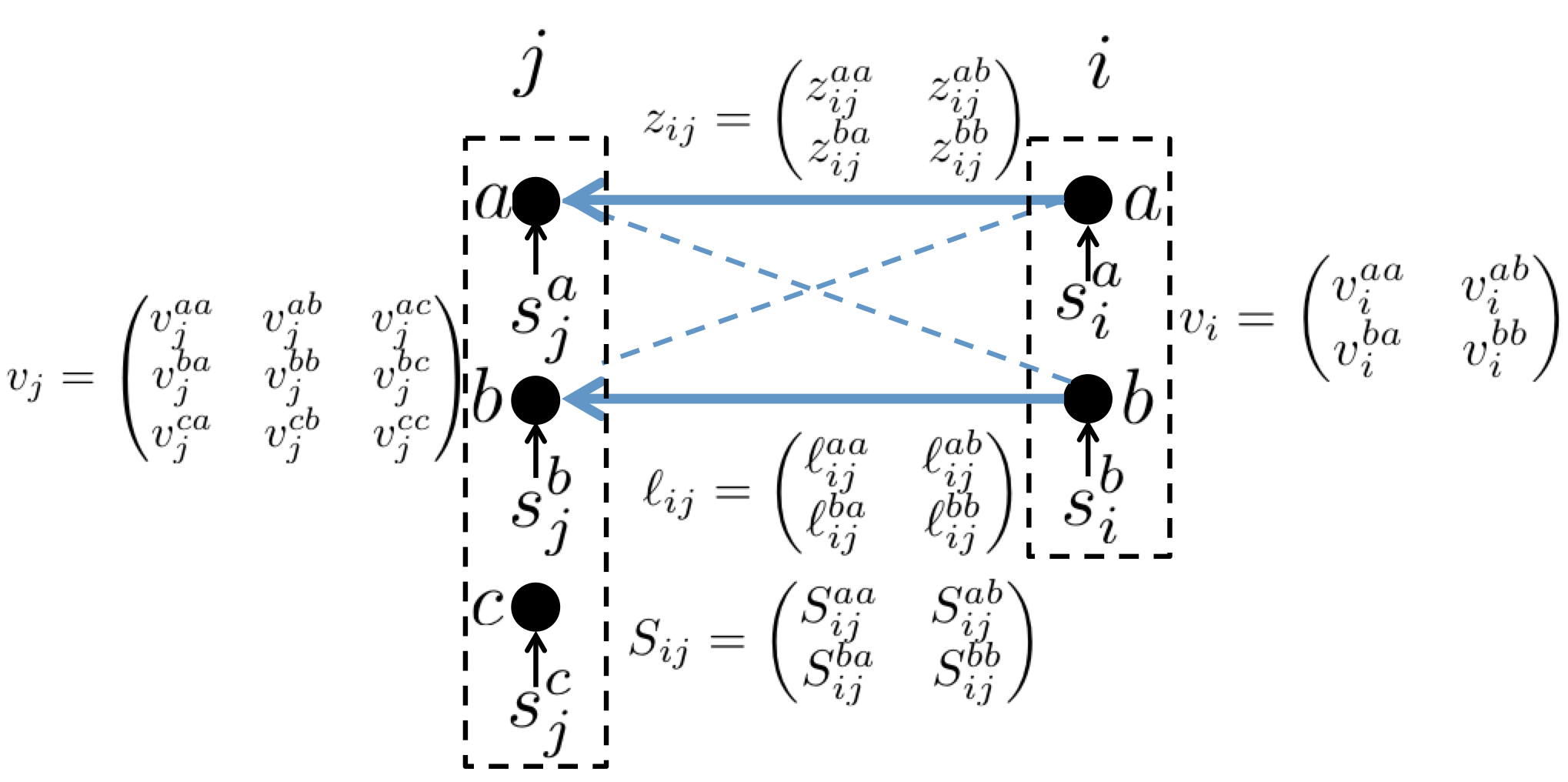}
\caption{Notations of a two bus network.}
\label{fig:notation}
\end{figure}

Let $a,b,c$ denote the three phases of the network. For each bus $i\in\mathN$, let $\Phi_i\subseteq\{a,b,c\}$ denote the set of phases. In a typical distribution network, the set of phases for bus $i$ is a subset of the phases of its parent and superset of the phases of its children, i.e. $\Phi_i\subseteq\Phi_{A_i}$ and $\Phi_j\subseteq \Phi_{i}$ for $j\in C_i$. On each phase $\phi\in\Phi_i$, let $V_i^{\phi}\in\mathbb{C}$ denote the complex voltage and $s_i^{\phi}:=p_i^{\phi}+jq_i^{\phi}$ denote the complex power injection. Denote $V_i:=(V_i^\phi,\phi\in\Phi_i)\in\mathbb{C}^{|\Phi_i|}$, $s_i:=(s_i^\phi,\phi\in\Phi_i)\in\mathbb{C}^{|\Phi_i|}$ and $v_i:=V_i^HV_i\in\mathbb{C}^{|\Phi_i|\times|\Phi_i|}$. For each line $i\in\mathN_+$ connecting bus $i$ and its ancestor $A_i$, the set of phases is $\Phi_i\cap\Phi_{A_i}=\Phi_i$ since $\Phi_i\subseteq\Phi_{A_i}$. On each phase $\phi\in\Phi_i$, let $I_i^{\phi}\in\mathbb{C}$ denote the complex branch current. Denote $I_i:=(I_i^\phi,\phi\in\Phi_i)\in\mathbb{C}^{|\Phi_i|}$, $\ell_i:=I_iI_i^H\in\mathbb{C}^{|\Phi_i|\times|\Phi_i|}$ and $S_i:=V_iI_i^H\in\mathbb{C}^{|\Phi_i|\times|\Phi_i|}$. 
Some notations are summarized in Fig. \ref{fig:notation}. A variable without a subscript denotes the set of variables with appropriate components, as summarized below.

\begin{center}
\begin{tabular}{|c|c|}
\hline
$v:=(v_i,i\in\mathN)$ & $s:=(s_i,i\in\mathN)$ \\
\hline
$\ell:=(\ell_i,i\in\mathN_+)$ & $S:=(S_i,i\in\mathN_+)$\\
\hline
\end{tabular}
\end{center}

Branch flow model is first proposed in \cite{Baran1989a,Baran1989b} for balanced radial networks. It has better numerical stability than bus injection model and has been advocated for the design and operation for radial distribution network, \cite{Farivar-2013-BFM-TPS, gan2014exact,li2012demand,peng2014distributed}. In \cite{gan2014convex}, it is first generalized to unbalanced radial networks and uses a set of variables $(v,s,\ell,S)$. Given a radial network $\mathT$, the branch flow model for unbalanced network is defined by:
\begin{subequations} \label{eq:bfm}
\begin{align}
& \mathP_i(v_{A_i})=v_i-z_iS_i^H-S_iz_i^H+z_i\ell_iz_i^H & i\in\mathN_+  \label{eq:bfm1}\\
& s_i=-\text{diag}\left(\sum_{i\in C_i}\mathP_i(S_j-z_j\ell_j)-S_i\right)  & i\in\mathN  \label{eq:bfm2}\\
& \begin{pmatrix}
v_i & S_i\\
S_i^H & \ell_i
\end{pmatrix}\in \mathbb{S}_+ &  i\in\mathN_+  \label{eq:bfm3}\\
& \text{rank}\begin{pmatrix}
v_i & S_i\\
S_i^H & \ell_i
\end{pmatrix}=1  &  i\in\mathN_+  \label{eq:bfm4} 
\end{align}
\end{subequations}
where $\mathP_i(v_{A_i})$ denote projecting $v_{A_i}$ to the set of phases on bus $i$ and $\mathP_i(S_j-z_j\ell_j)$ denote lifting the result of $S_j-z_j\ell_j$ to the set of phases $\Phi_i$ and filling the missing phase with $0$, e.g. if $\Phi_{A_i}=\{a,b,c\}$, $\Phi_i=\{a,b\}$ and $\Phi_j=\{a\}$, then 
\bqn
\mathP_i(v_{A_i})&:=&\begin{pmatrix}
v_{A_i}^{aa} & v_{A_i}^{ab} \\
v_{A_i}^{ba} & v_{A_i}^{bb}
\end{pmatrix} \\
\mathP_i(S_j-z_j\ell_j)&:=&\begin{pmatrix}
S_{j}^{aa}-z_j^{aa}\ell_j^{aa} & 0 \\
0 & 0
\end{pmatrix} 
\eqn

Given a vector $(v,s,\ell,S)$ that satisfies \eqref{eq:bfm}, it is proved in \cite{gan2014convex} that the bus voltages $V_i$ and branch currents $I_i$ can be uniquely determined if the network is a tree. Hence this model  \eqref{eq:bfm} is equivalent to a full unbalanced AC power flow model. See \cite[Section IV]{gan2014convex} for details.

\subsection{OPF and SDP relaxation}\label{sec:opfsdp}

The OPF problem seeks to optimize certain objective, e.g. total power loss or generation cost, subject to unbalanced power flow equations \eqref{eq:bfm} and various operational constraints. 
We consider an objective function of the following form: 
\bq\label{eq:objective}
F(s):=\sum_{i\in\mathN}f_i(s_i):=\sum_{i\in\mathN}\sum_{\phi\in\Phi_i}f_i^\phi(s_i^\phi).
\eq
For instance,
\bi
\item to minimize total line loss, we can set for each $\phi\in\Phi_i$, $i\in\mathN$,
\bq\label{eq:opf::unbalance::I2}
f_i^\phi(s_i^\phi) =  p_i^\phi.
\eq
\item to minimize generation cost, we can set for each $i\in\mathN$,
\bq\label{eq:opf::unbalance::I1}
f_i^\phi(s_i^\phi)=(\frac{\alpha_i^\phi}{2} (p_i^\phi)^2+\beta_i^\phi p_i^\phi),
\eq
where $\alpha_i^\phi,\beta_i^\phi>0$ depend on the load type on bus $i$, e.g. $\alpha_i^\phi=0$ and $\beta_i^\phi=0$ for bus $i$ where there is no generator and for generator bus $i$, the corresponding $\alpha_i^\phi,\beta_i^\phi$ depends on the characteristic of the generator. 
\ei

For each bus $i\in\mathN$, there are two operational constraints on each phase $\phi\in\Phi_i$. First, the power injection $s_i^\phi$ is constrained to be in a injection region $\mathI_i^\phi$, i.e.
\bq\label{eq:operation1}
s_i^\phi\in\mathI_i^\phi \ \ \text{for } \phi\in\Phi_i \text{ and }i\in\mathN
\eq

The feasible power injection region $\mathI_i^\phi$ is determined by the controllable loads attached to phase $\phi$ on bus $i$. Some common controllable loads are:

\bi
\item For controllable load, whose real power can vary within $[\underline p_i,\overline p_i]$ and reactive power can vary within $[\underline q_i,\overline q_i]$, the injection region $\mathI_i$ is 
\begin{subequations}
\bq\label{eq:S2}
\mathI_i^\phi=\{p+\mathbf{i}q\mid p\in[ \underline p_i,  \overline p_i], q\in[\underline q_i,\overline q_i] \}\subseteq \mathbb{C}
\eq
For instance, the power injection of each phase $\phi$ on substation bus $0$ is unconstrained, thus $\underline p_i,\underline q_i=-\infty$ and $ \overline p_i, \overline q_i=\infty$. 
\item For solar panel connecting the grid through a inverter with nameplate $\overline s_i^\phi$, the injection region $\mathI_i$ is 
\bq\label{eq:S1}
\mathI_i^\phi=\{p+\mathbf{i}q\mid p\geq 0, p^2+q^2\leq  (\overline s_i^\phi)^2\}\subseteq\mathbb{C}
\eq
\end{subequations}
\ei

Second, the voltage magnitude needs to be maintained within a prescribed region. Note that the diagonal element of $v_i$ describes the voltage magnitude square on each phase $\phi\in\Phi_i$. Thus the constraints can be written as 
\bq\label{eq:operation2}
\underline v_i^\phi \leq v_i^{\phi\phi}\leq \overline v_i^\phi \ \ i\in\mathN,
\eq
where $v_i^{\phi\phi}$ denotes the $\phi_{th}$ diagonal element of $v_i$. Typically the voltage magnitude at substation buses is assumed to be fixed at a prescribed value,  i.e. $\underline v_0^{\phi}=\overline v_0^{\phi}$ for $\phi\in\Phi_0$. At other load buses $i\in\mathN_+$, the voltage magnitude is typically allowed to deviate by $5\%$ from its nominal value, i.e. $\underline v_i^{\phi}=0.95^2$ and $\overline v_i^{\phi}=1.05^2$ for $\phi\in\Phi_i$. 

To summarize, the OPF problem for unbalanced radial distribution networks is: 
\bq
\text{\bf OPF: }\min && \sum_{i\in\mathN}\sum_{\phi\in\Phi_i}f_i^\phi(s_i^\phi)\nonumber\\
\mathrm{over} && v,s,S,\ell \label{eq:opf}\\
\mathrm{s.t.} &&  \eqref{eq:bfm} \text{ and } \eqref{eq:operation1}-\eqref{eq:operation2}\nonumber
\eq

The OPF problem \eqref{eq:opf} is nonconvex due to the rank constraint \eqref{eq:bfm4}. In \cite{gan2014convex}, an SDP relaxation for \eqref{eq:opf} is obtained by removing the rank constraint \eqref{eq:bfm4}, resulting in a semidefinite program (SDP): 
\bq
\text{\bf ROPF: }\min && \sum_{i\in\mathN}\sum_{\phi\in\Phi_i}f_i^\phi(s_i^\phi)\nonumber\\
\mathrm{over} && v,s,S,\ell \label{eq:ropf}\\
\mathrm{s.t.} &&  \eqref{eq:bfm1}-\eqref{eq:bfm3}  \text{ and }  \eqref{eq:operation1}-\eqref{eq:operation2}\nonumber
\eq

Clearly the relaxation ROPF \eqref{eq:ropf} provides a lower bound for the original OPF problem \eqref{eq:opf} since the original feasible set is enlarged. The relaxation is called \emph{exact} if every optimal solution of ROPF satisfies the rank constraint \eqref{eq:bfm4} and hence is also optimal for the original OPF problem. It is shown empirically in \cite{gan2014convex} that the relaxation is exact for all the tested distribution networks, including IEEE test networks \cite{kersting1991radial} and some real distribution feeders.

\section{Distributed Algorithm}\label{sec:alg}

We assume SDP relaxation is exact and  develop in this section
a distributed algorithm that solves the ROPF problem. 
We first design a distributed algorithm for a broad class of optimization problem through alternating direction method of multipliers (ADMM). We then apply the proposed algorithm on the ROPF problem, and show that the optimization subproblems can be solved efficiently either through closed form solutions or eigen-decomposition of a $6\times 6$ matrix.

\subsection{Preliminary: ADMM}

ADMM blends the decomposability of dual decomposition with the superior convergence properties of the method of multipliers
\cite{boyd2011distributed}. It solves optimization problem of the form\footnote{This is a special case with simpler constraints of the general form introduced in \cite{boyd2011distributed}. The $z$ variable used in \cite{boyd2011distributed} is replaced by $y$ since $z$ represents impedance in power systems.}:
\bq
\min_{x,y} && f(x)+g(y) \nonumber \\
\text{s.t.} && x\in\mathcal{K}_x, \ \ y\in\mathcal{K}_y \label{eq:admm}\\
&& x=y \nonumber
\eq
where $f(x),g(y)$ are convex functions and $\mathcal{K}_x,\mathcal{K}_y$ are convex sets. Let $\lambda$ denote the Lagrange
 multiplier for the constraint $x=y$. Then the augmented Lagrangian is defined as
\bq\label{eq:agumentlag}
L_\rho(x,y,\lambda):=f(x)+g(y)+\langle \lambda, x-y\rangle+\frac{\rho}{2}\|x-y\|_2^2,
\eq
where  $\rho\geq 0$ is a constant. When $\rho=0$, the augmented Lagrangian degenerates to 
the standard Lagrangian. At each iteration $k$, ADMM consists of the iterations:
\begin{subequations}\label{eq:update}
\bq
x^{k+1}&\in&\arg\min_{x\in\mathcal{K}_x} L_\rho(x,y^{k},\lambda^k)\label{eq:xupdate}\\
y^{k+1}&\in&\arg\min_{y\in\mathcal{K}_y} L_\rho(x^{k+1},y,\lambda^k)\label{eq:zupdate}\\
\lambda^{k+1}&=&\lambda^{k}+\rho(x^{k+1}-y^{k+1}).\label{eq:mupdate}
\eq
\end{subequations}
Specifically, at each iteration, ADMM first updates $x$ based on \eqref{eq:xupdate}, then updates $y$ based on \eqref{eq:zupdate}, and after that updates the multiplier $\lambda$ based on \eqref{eq:mupdate}. Compared to dual decomposition, ADMM is guaranteed to converge to an optimal solution under less restrictive conditions. Let
\begin{subequations}\label{eq:feasible}
\bq
r^k&:=&\|x^{k}-y^{k}\|_2 \label{eq:pfeasible} \\ 
s^k&:=&\rho\|y^{k}-y^{k-1}\|_2, \label{eq:dfeasible}
\eq
\end{subequations}
which can be viewed as the residuals for primal and dual feasibility, respectively. They converge to $0$ at optimality and are usually used as metrics of convergence in the experiment. Interested readers may refer to \cite[Chapter 3]{boyd2011distributed} for details. 

In this paper, we generalize the above standard ADMM \cite{boyd2011distributed} such that the optimization subproblems can be solved efficiently for our ROPF problem. Instead of using the quadratic penalty term $\frac{\rho}{2}\|x-y\|_2^2$ in \eqref{eq:agumentlag}, we will use a more general quadratic penalty term: $\|x-y\|_{\Lambda}^2$, where $\|x-y\|_{\Lambda}^2:=(x-y)^H\Lambda(x-y)$ and $\Lambda$ is a positive diagonal matrix. Then the augmented Lagrangian becomes
\bq\label{eq:admm::augmentlag}
L_{\rho}(x,y,\lambda):=f(x)+g(y)+\langle \lambda, x-y\rangle+\frac{\rho}{2}\|x-y\|_{\Lambda}^2.
\eq
The convergence result in  \cite[Chapter 3]{boyd2011distributed} carries over directly to this general case. 

\subsection{ADMM based Distributed Algorithm}\label{sec:admma_dalg}

In this section, we will design an ADMM based distributed algorithm for a broad class of optimization problem, of which the ROPF problem is a special case. Consider the following optimization problem:
\begin{subequations}\label{eq:admm::opt1}
\bq
\min && \sum_{i\in\mathN} f_i(x_i) \label{eq:admm::opt1::obj}\\
\mathrm{over} && \{x_i\mid i\in\mathN\} \\
\mathrm{s.t.}  && \sum_{j\in N_i}A_{ij}x_j = 0 \quad\mathrm{for} \quad i\in\mathN \label{eq:admm::couple}\\
&& x_i\in \cap_{r=0}^{R_i}\mathK_{ir} \quad\mathrm{for} \quad i\in\mathN, \label{eq:admm::local}
\eq
\end{subequations}
where for each $i\in\mathN$, $x_i$ is a complex vector, $f_i(x_i)$ is a convex function, $\mathK_{ir}$ is a convex set, and $A_{ij}$ $(j\in N_i, i\in\mathN)$ are matrices with appropriate dimensions. A broad class of graphical optimization problems (including ROPF) can be formulated as \eqref{eq:admm::opt1}. Specifically, each node $i\in\mathN$ is associated with some local variables stacked as $x_i$, which belongs to an intersection of $R_i+1$ local feasible sets $\mathK_{ir}$ and has a cost objective function $f_i(x_i)$. Variables in node $i$ are coupled with variables from their neighbor nodes in $N_i$ through linear constraints \eqref{eq:admm::couple}. The objective then is to solve a minimal total cost across all the nodes.

The goal is to develop a distributed algorithm that solves \eqref{eq:admm::opt1} such that each node $i$ solve its own subproblem and only exchange information with its neighbor nodes $N_i$. In order to transform \eqref{eq:admm::opt1} into the form of standard ADMM \eqref{eq:admm}, we need to have two sets of variables $x$ and $y$. We introduce two sets of slack variables as below:
\bee
\item $x_{ir}$. It represents a copy of the original variable $x_i$ for $1\leq r \leq R_i$. For convenience, denote the original $x_i$ by $x_{i0}$.
\item $y_{ij}$. It represents the variables in node $i$ observed at node $j$, for $j\in N_i$.
\eee
Then \eqref{eq:admm::opt1} can be reformulated as
\begin{subequations}\label{eq:admm::opt4}
\bq
\min && \sum_{i\in\mathN} f_i(x_{i0}) \label{eq:admm::opt4::obj}\\
\mathrm{over} && x=\{x_{ir}\mid 0\leq r\leq R_i,i\in\mathN\} \nonumber\\ 
&&y=\{y_{ij}\mid j\in N_i,i\in\mathN\} \nonumber  \\
\mathrm{s.t.}  && \sum_{j\in N_i}A_{ij}y_{ji} = 0 \quad\mathrm{for} \quad i\in\mathN \label{eq:admm::opt4::couple}\\
&& x_{ir}\in \mathK_{ir} \quad\mathrm{for} \quad   0\leq r\leq R_i \ \  i\in\mathN \label{eq:admm::opt4::local}\\
&& x_{ir}=y_{ii} \quad\mathrm{for} \quad   1\leq r\leq R_i \ \  i\in\mathN \label{eq:admm::opt4::consensus1} \\
&& x_{i0}=y_{ij} \quad\mathrm{for} \quad   j\in N_i \ \  i\in\mathN , \label{eq:admm::opt4::consensus2} 
\eq
\end{subequations}
where $x$ and $y$ represent the two groups of variables in standard ADMM. Note that the consensus constraints \eqref{eq:admm::opt4::consensus1} and \eqref{eq:admm::opt4::consensus2} force all the duplicates $x_{ir}$ and $y_{ij}$ are the same. Thus its solution $x_{i0}$ is also optimal to the original problem \eqref{eq:admm::opt1}. \eqref{eq:admm::opt4} falls into the general ADMM form \eqref{eq:admm}, where \eqref{eq:admm::opt4::couple} corresponds to $\mathK_y$, \eqref{eq:admm::opt4::local} corresponds to $\mathK_x$, and \eqref{eq:admm::opt4::consensus1} and \eqref{eq:admm::opt4::consensus2} are the consensus constraints that relates $x$ and $y$.

Following the ADMM procedure, we relax the consensus constraints \eqref{eq:admm::opt4::consensus1} and \eqref{eq:admm::opt4::consensus2}, whose Lagrangian multipliers are denoted by $\lambda_{ir}$ and $\mu_{ij}$, respectively. The generalized augmented Lagrangian then can be written as 
\begin{align}\label{eq:admm::opt4::augment}
&L_{\rho}(x,y,\lambda,\mu)\\
=&\sum_{i\in\mathN} \left(\sum_{r=1}^{R_i}\left(\langle\lambda_{ir}, x_{ir}-y_{ii}\rangle+ \frac{\rho}{2}\|x_{ir}-y_{ii}\|_{\Lambda_{ir}}^2\right)+\right. \nonumber\\
&\left.f_i(x_{i0})+\sum_{j\in N_i}\left(\langle\mu_{ij},x_{i0}-y_{ij}\rangle +\frac{\rho}{2}\|x_{i0}-y_{ij}\|_{M_{ij}}^2\right) \right).\nonumber
\end{align}
where the parameter $\Lambda_{ir}$ and $M_{ij}$ depend on the problem we will show how to design them in section \ref{sec:dalg}.

Next, we show that both the $x$-update \eqref{eq:xupdate} and $y$-update \eqref{eq:zupdate} can be solved in a distributed manner, i.e. both of them can be decomposed into local subproblems that can be solved in parallel by each node $i$ with only neighborhood communications.

First, we define the set of local variables for each node $i$, denoted by $\mathcal{A}_i$, which includes its own duplicates $x_{ir}$ and the associated multiplier $\lambda_{ir}$ for $0\leq r\leq R_i$, and the ``observations'' $y_{ji}$ of variables from its neighbor $N_i$ and the associated multiplier $\mu_{ji}$, i.e.
\bq\label{eq:localvar}
\mathcal{A}_i:=\{x_{ir},\lambda_{ir}\mid 0\leq r\leq R_i\}\cup \{y_{ji},\mu_{ji}\mid j\in N_i\}.
\eq
Next, we show how does each node $i$ update $\{x_{ir}\mid 0\leq r\leq R_i\}$ in the $x$-update and $\{y_{ji}\mid j\in N_i\}$ in the $y$-update.

\begin{figure}
\centering
\subfloat[$x$-update]{
	\includegraphics[scale=0.2]{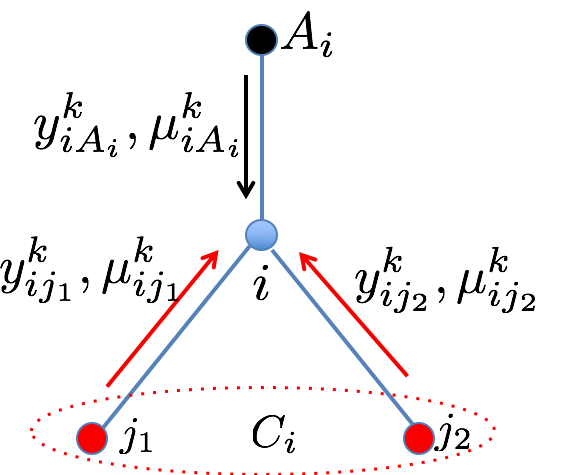}
	\label{fig:msg_x}
	}
\subfloat[$y$-update] {
	\includegraphics[scale=0.2]{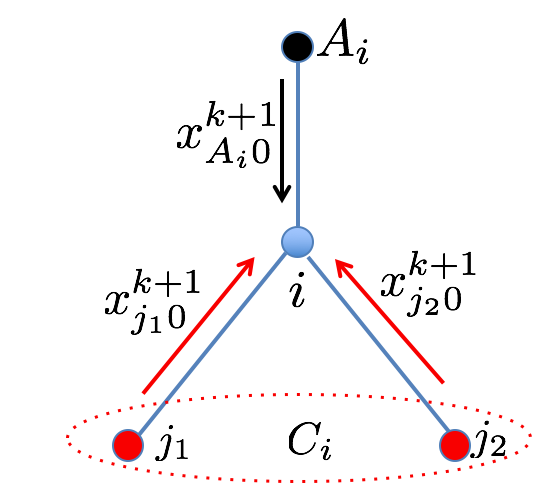}
	\label{fig:msg_z}
	}
\caption{Message exchanges in the $x$ and $y$-update for node $i$.}
\label{fig:msg_update}
\end{figure}

In the $x$-update at each iteration $k$, the optimization subproblem that updates $x^{k+1}$ is 
\bq\label{eq:admm:opt4:xupdate}
\min_{x\in\mathK_x} L_\rho(x, y^k, \lambda^k,\mu^k),
\eq
where the constraint $\mathcal{K}_x$ is the Cartesian product of $\mathK_{ir}$, i.e.
\bqn
\mathK_x:=\otimes_{i\in\mathN}\otimes_{r=0}^{R_i}\mathK_{ir}.
\eqn
The objective can be written as a sum of local objectives as shown below
{\small
\begin{align*}
&L_{\rho}(x,y^k,\lambda^k,\mu^k)\\
=&\sum_{i\in\mathN} \left(\sum_{r=1}^{R_i}\left(\langle\lambda_{ir}^{k}, x_{ir}-y_{ii}^{k}\rangle+ \frac{\rho}{2}\|x_{ir}-y_{ii}^{k}\|_{\Lambda_{ir}}^2\right)+\right. \nonumber\\
&\left.f_i(x_{i0})+\sum_{j\in N_i}\left(\langle\mu_{ij}^{k},x_{i0}-y_{ij}^{k}\rangle +\frac{\rho}{2}\|x_{i0}-y_{ij}^{k}\|_{M_{ij}}^2\right) \right)\nonumber\\
=& \sum_{i\in\mathN}\sum_{r=0}^{R_i}H_{ir}(x_{ir}) - \sum_{i\in\mathN}\left(\sum_{r=0}^{R_i}\langle \lambda_{ir}^k,y_{ii}^k\rangle+\sum_{j\in N_i}\langle \mu_{ij}^k,y_{ij}^k\rangle\right),
\end{align*}
}
where the last term is independent of $x$ and 
\begin{align} \label{eq:admm:opt4:hix}
&H_{ir}(x_{ir})=\\
&\begin{cases}
   f_i(x_{i0})+\sum_{j\in N_i}\left(\langle\mu_{ij}^{k},x_{i0}\rangle +\frac{\rho}{2}\|x_{i0}-y_{ij}^{k}\|_{M_{ij}}^2\right)& r=0\\
         \langle\lambda_{ir}^{k},x_{ir}\rangle+ \frac{\rho}{2}\|x_{ir}-y_{ii}^{k}\|_{\Lambda_{ir}}^2        & r > 0
\end{cases}.\nonumber
\end{align}
Then the problem \eqref{eq:admm:opt4:xupdate} in the $x$-update can be written explicitly as
\bq
\min && \sum_{i\in\mathN}\sum_{r=0}^{R_i} H_{ir}(x_{ir}) \nonumber\\
\mathrm{over} && x=\{x_{ir}\mid 0\leq r\leq R_i, i\in\mathN\}  \label{eq:admm::opt2::xupdate}\\
\mathrm{s.t.} && x_{ir}\in\mathK_{ir} \quad\mathrm{for} \quad 0\leq r\leq R_i, \ i\in\mathN, \nonumber
\eq
where both the objective and constraint are separable for $0\leq r\leq R_i$ and $i\in\mathN$. Thus it can be decomposed into $\sum_{i\in\mathN}(R_i+1)$ independent problems that can be solved in parallel. There are $R_i+1$ problems associated with each node $i$ and the $r_{th}$ $(0\leq r\leq R_i)$ one can be simply written as 
\bq\label{eq:xupdatenode}
\min_{x_{ir}\in\mathK_{ir}} \ \ H_{ir}(x_{ir})
\eq
whose solution is the new update of variables $x_{ir}$ for node $i$. In the above problem, the constants $y_{ij}^k,\mu_{ij}^k\in\mathcal{A}_j$ are not local to $i$ and stored in $i$'s neighbors $j\in N_i$. Therefore, each node $i$ needs to collect $(y_{ij},\mu_{ij})$ from all of its neighbors prior to solving \eqref{eq:xupdatenode}. The message exchanges is illustrated in Figure \ref{fig:msg_x}.

In the $y$-update, the optimization problem that updates $y^{k+1}$ is 
\bq\label{eq:admm:opt4:yupdate1}
\min_{y\in\mathK_y} L_\rho(x^{k+1}, y, \lambda^k,\mu^k)
\eq
where the constraint set $\mathK_y$ can be represented as a Cartesian product of $|\mathN|$ disjoint sets, i.e.
\bqn
\mathK_y:=\otimes_{i\in\mathN}\{y_{ji}, j\in N_i\mid  \sum_{j\in N_i}A_{ij}y_{ji} = 0 \}.
\eqn

The objective can be written as a sum of local objectives as below.
{\small
\begin{align*}
&L_{\rho}(x^{k+1},y,\lambda^k,\mu^k)\\
=&\sum_{i\in\mathN} \left(\sum_{r=1}^{R_i}\left(\langle\lambda_{ir}^{k}, x_{ir}^{k+1}-y_{ii}\rangle+ \frac{\rho}{2}\|x_{ir}^{k+1}-y_{ii}\|_{\Lambda_{ir}}^2\right)+\right. \nonumber\\
&\left.f_i(x_{i0}^{k+1})+\sum_{j\in N_i}\left(\langle\mu_{ji}^{k},x_{j0}^{k+1}-y_{ji}\rangle +\frac{\rho}{2}\|x_{j0}^{k+1}-y_{ji}\|_{M_{ji}}^2\right) \right)\nonumber\\
=&\sum_{i\in\mathN}G_i(\{y_{ji}\mid j\in N_i\}) + \\
&\sum_{i\in\mathN}\left(f_i(x_{i0}^{k+1})+\sum_{r=0}^{R_i}\langle \lambda_{ir}^{k}, x_{ir}^{k+1}\rangle+\sum_{j\in N_i}\langle \mu_{ji}^{k}, x_{j0}^{k+1}\rangle\right),
\end{align*}
}
where the last term is independent of $y$ and 
\begin{align*}
&G_i(\{y_{ji}\mid j\in N_i\})=\sum_{r=1}^{R_i}\left(-\langle\lambda_{ir}^{k},y_{ii}\rangle+ \frac{\rho}{2}\|x_{ir}^{k+1}-y_{ii}\|_{\Lambda_{ir}}^2\right)\\
&\qquad +\sum_{j\in N_i}\left(-\langle\mu_{ji}^{k},y_{ji}\rangle +\frac{\rho}{2}\|x_{j0}^{k+1}-y_{ji}\|_{M_{ji}}^2\right).
\end{align*}
Then the problem \eqref{eq:admm:opt4:yupdate1} in the $y$-update can be written explicitly as
\bqn
\min && \sum_{i\in\mathN}G_i(\{y_{ji}\mid j\in N_i\})\\
\mathrm{over} && y=\{\{y_{ji}\mid j\in N_i\}\mid i\in\mathN\}\\
\mathrm{s.t.} && \sum_{j\in N_i}A_{ij}y_{ji} = 0, \ \ i \in\mathN
\eqn
which can be decomposed into $|\mathN|$ subproblems and the subproblem for node $i$ is 
\bq
\min && G_i(\{y_{ji}\mid j\in N_i\})\nonumber\\
\mathrm{over} && \{y_{ji}\mid j\in N_i\}  \label{eq:admm::opt4::zupdate1}\\
\mathrm{s.t.} && \sum_{j\in N_i}A_{ij}y_{ji} = 0, \nonumber
\eq
whose solution is the new update of $\{y_{ji}\mid j\in N_i\}\in\mathcal{A}_i$. In \eqref{eq:admm::opt4::zupdate1}, the constants
 $x_{j0}\in\mathcal{A}_j$ are stored in $i$'s neighbor $j\in N_i$. Hence, each node $i$ needs to collect $x_{j0}$ from all of its neighbor prior to solving \eqref{eq:admm::opt4::zupdate1}. The message exchanges in the $y$-update is illustrated in Figure \ref{fig:msg_z}.

The problem \eqref{eq:admm::opt4::zupdate1} can be solved with closed form solution. we stack the real and imaginary part of the variables $\{y_{ji}\mid j\in N_i\}$ in a vector with appropriate dimensions and denote it as $\tilde y$. Then \eqref{eq:admm::opt4::zupdate1} takes the following form:
\bq
\min && \frac{1}{2}{\tilde y}^TM {\tilde y}+c^T{\tilde y} \nonumber\\
\mathrm{over} && \tilde y  \label{eq:admm::opt2::zupdate1} \\
\mathrm{s.t.} && \tilde A{\tilde y}=0,  \nonumber
\eq
where $M$ is a positive diagonal matrix, $\tilde A$ is a full row rank real matrix, and $c$ is a real vector. $M,c,A$ are derived from \eqref{eq:admm::opt4::zupdate1}. There exists a closed form expression for \eqref{eq:admm::opt2::zupdate1} given by
\bq \label{eq:admm::opt2::zsol}
\tilde y=\left(M^{-1}\tilde A^T(\tilde AM^{-1}\tilde A^T)^{-1}\tilde AM^{-1}-M^{-1}\right)c.
\eq

In summary, the original problem \eqref{eq:admm::opt1} is decomposed into local subproblems that can be solved in a distributed manner using ADMM. At each iteration, each node $i$ solves \eqref{eq:xupdatenode} in the $x$-update and \eqref{eq:admm::opt2::zupdate1} in the $y$-update. There exists a closed form solution to the subproblem \eqref{eq:admm::opt2::zupdate1} in the $y$-update as shown in \eqref{eq:admm::opt2::zsol}, and hence whether the original problem \eqref{eq:admm::opt1} can be solved efficiently in a distributed manner depends on the existence of efficient solutions to the subproblems \eqref{eq:xupdatenode} in the $x$-update, which depends on the realization of both the objectives $f_i(x_i)$ and the constraint sets $\mathK_{ir}$.

Next, we show the ROPF problem \eqref{eq:ropf} is a special case of \eqref{eq:admm::opt1}, hence can be solved in a distributed manner using the above method. In particular, we show the corresponding subproblems in the $x$-update can be solved efficiently.

\subsection{Application on OPF problem}\label{sec:dalg}
We assume the SDP relaxation is exact and now derive a distributed algorithm for solving ROPF \eqref{eq:ropf}. Using the ADMM based algorithm developed in Section \ref{sec:admma_dalg}, the global ROPF problem is decomposed into local subproblems that can be solved in a distributed manner with only neighborhood communication. Note that the subproblems in the $y$-update for each node $i$ can always been solved with closed form solution, we only need to develop an efficient solution for the subproblems \eqref{eq:xupdatenode} in the $x$-update for the ROPF problem. In particular, we provide a sufficient condition, which holds in practice, for the existence of efficient solutions to all the optimization subproblems. Compared with existing methods, e.g. \cite{devane2013stability,li2012demand,dall2013distributed,kraning2013dynamic,sun2013fully}, that use generic iterative optimization solver to solve each subproblem, the computation time is improved by more than 100 times.

The ROPF problem defined in \eqref{eq:ropf} can be written explicitly as
\begin{subequations}\label{eq:ropfexplicit}
\begin{align}
\min & \ \ \sum_{i\in \mathN} \sum_{\phi\in\Phi_i}f_i^\phi(s_i^\phi) \\
\mathrm{over} & \ \ v,s,S,\ell \nonumber \\
\mathrm{s.t.}  & \ \ \mathP_i(v_{A_i})=v_i-z_iS_i^H-S_iz_i^H+z_i\ell_iz_i^H  \ i\in\mathN  \label{eq:bfm21}\\
& s_i=-\text{diag}\left(\sum_{i\in C_i}\mathP_i(S_j-z_j\ell_j)-S_i\right)  \quad i\in\mathN  \label{eq:bfm22}\\
& \begin{pmatrix}
v_i & S_i\\
S_i^H & \ell_i
\end{pmatrix}\in \mathbb{S}_+  \qquad\qquad\qquad\qquad\quad \ \  i\in\mathN  \label{eq:bfm23}\\
& \ \ s_i^\phi\in\mathI_i^\phi \qquad\qquad\qquad\qquad\qquad \phi\in\Phi_i, \ i\in\mathN   \label{eq:bfm24}\\
& \ \ \underline v_i^\phi \leq v_i^{\phi\phi}\leq \overline v_i^{\phi}  \qquad\qquad\quad\qquad \phi\in\Phi_i, \ i\in\mathN   \label{eq:bfm25}
\end{align}
\end{subequations}
Denote 
\begin{align}
x_i:=&\{v_i,s_i,S_i,\ell_i\} \\
\mathK_{i0}:=&\{x_i\mid 
\begin{pmatrix}
v_i & S_i\\
S_i^H & \ell_i
\end{pmatrix}\in \mathbb{S}_+, \{s_i^\phi\in\mathcal{I}_i^\phi\mid \phi\in \Phi_i\}\} \label{eq:distOPFu::ki0}\\
\mathK_{i1}:=&\{x_i\mid \underline v_i^\phi \leq v_i^{\phi\phi}\leq \overline v_i^{\phi}, \phi\in \Phi_i\} \label{eq:distOPFu::ki1}
\end{align}
Then \eqref{eq:ropfexplicit} takes the form of \eqref{eq:admm::opt1} with $R_i=1$ for all $i\in\mathN$, where \eqref{eq:bfm21}--\eqref{eq:bfm22} correspond to \eqref{eq:admm::couple} and \eqref{eq:bfm23}--\eqref{eq:bfm25} correspond to \eqref{eq:admm::local}. Then we have the following theorem, which provides a sufficient condition for the existence of an efficient solution to \eqref{eq:xupdatenode}.
\begin{theorem}\label{thm:distOPFu:closedformsol}
Suppose there exists a closed form solution to the following optimization problem for all $i\in\mathN$ and $\phi\in\Phi_i$
\bq
\min && f_i^\phi\left(s\right)+\frac{\rho}{2}\left\|s^\phi-\hat s^\phi\right\|_2^2 \nonumber\\
\mathrm{over} && s \in \mathI_i^\phi   \label{eq:distOPFu:closedformsol} 
\eq
given any constant $\hat s^\phi$ and $\rho$, then the subproblems for ROPF in the $x$-update \eqref{eq:xupdatenode} can be solved via either closed form solutions or eigen-decomposition of a $6\times 6$ hermitian matrix.
\end{theorem}
\begin{IEEEproof}
We will prove Theorem \ref{thm:distOPFu:closedformsol} through elaborating the procedure to solve \eqref{eq:xupdatenode}.
\end{IEEEproof}

Recall that there is always a closed form solution to the optimization subproblem \eqref{eq:admm::opt4::zupdate1} in the $y$-update, if the objective function $f_i^\phi\left(s^\phi\right)$ and injection region $\mathI_i^\phi$ satisfy the sufficient condition in Theorem \ref{thm:distOPFu:closedformsol}, all the subproblems can be solved efficiently.

\begin{remark}
In practice, the objective function $f_i^\phi(s)$, usually takes the form of $f_i^\phi\left(s\right):=\frac{\alpha_i}{2} p^2+\beta_i p$, which models both line loss and generation cost minimization as discussed in Section \ref{sec:opfsdp}. For the injection region $\mathI_i^\phi$, it usually takes either \eqref{eq:opf::unbalance::I2} or \eqref{eq:opf::unbalance::I1}. It is shown in Appendix \ref{app:distOPFb::solver2} that there exist closed form solution for all of those cases. Thus \eqref{eq:distOPFu:closedformsol} can be solved efficiently for practical applications. 
\end{remark}

Following the procedure in Section \ref{sec:admma_dalg}, we introduce two set of slack variables: $x_{ir}$ and $y_{ij}$. Then the counterpart of \eqref{eq:admm::opt4} is 
{\scriptsize
\begin{subequations}\label{eq:distOPFu::eropf}
\begin{align}
\min & \ \ \sum_{i\in\mathN}\sum_{\phi\in\Phi_i}f_i^\phi((s_{i0}^\phi)^{(x)}) \\
\mathrm{over} & \ \ x:=\{x_{ir}\mid 0\leq r\leq 1, \ i\in\mathN \}  \nonumber\\
& \ \ y:=\{y_{ji} \mid j\in N_i, \ i\in\mathN\} \nonumber \\
\mathrm{s.t.}  & \ \ \mathP_i(v^{(y)}_{A_ii})=v^{(y)}_{ii}-z_i(S^{(y)}_{ii})^H-S^{(y)}_{ii}z_i^H+z_i\ell^{(y)}_{ii}z_i^H \  i\in\mathN \\
& \ \ s^{(y)}_{ii}=-\text{diag}\left(\sum_{i\in C_i}\mathP_i(S^{(y)}_{ji}-z_j\ell^{(y)}_{ji})-S^{(y)}_{ii}\right)  \qquad i\in\mathN \\
&\ \ \begin{pmatrix}
v_{i0}^{(x)} & S_{i0}^{(x)}\\
(S_{i0}^{(x)})^H & \ell_{i0}^{(x)}
\end{pmatrix}\in \mathbb{S}_+  \qquad \qquad \qquad \qquad \qquad \qquad \quad i\in\mathN  \\
&\ \  (s_{i0}^\phi)^{(x)}\in \mathI_i^\phi  \qquad\qquad\qquad \qquad \qquad \qquad \qquad \phi\in\Phi_i \text{ and }  i\in\mathN \\
&\ \ \underline v_i^\phi \leq (v_{i1}^{\phi\phi})^{(x)}\leq \overline v_i^\phi \quad\quad \qquad \qquad \qquad \qquad \phi\in\Phi_i  \text{ and }i\in\mathN  \\
& \ \ x_{ir} - y_{ii} = 0 \qquad\qquad\qquad \qquad \qquad \qquad \qquad r=1 \text{ and } i\in \mathN \label{eq:distOPFu:eropf::con1} \\
& \ \  x_{i0}-y_{ij}=0 \qquad\qquad\qquad \qquad \qquad \qquad \qquad j\in N_i \text{ and } i\in\mathN , \label{eq:distOPFu:eropf::con2}
\end{align}
\end{subequations}
}
where we put superscript $(\cdot)^{(x)}$ and $(\cdot)^{(y)}$ on each variable to denote whether the variable is updated in the $x$-update or $y$-update step. Note that each node $i$ does not need full information of its neighbor. Specifically, for each node $i$, only voltage information $v_{A_ii}^{(y)}$ is needed from its parent $A_i$ and branch power $S_{ji}^{(y)}$ and current $\ell_{ji}^{(y)}$ information from its children $j\in C_i$ based on \eqref{eq:distOPFu::eropf}. Thus, $y_{ij}$ contains only partial information about $x_{i0}$, i.e.
\bqn
y_{ij}&:=&\begin{cases}
(S_{ii}^{(y)}, \ell_{ii}^{(y)}, v_{ii}^{(y)}, s_{ii}^{(y)}) & j=i\\
(S_{iA_i}^{(y)}, \ell_{iA_i}^{(y)}) & j=A_i\\
(v_{ij}^{(y)}) & j\in C_i
\end{cases}.
\eqn
On the other hand, only $x_{i0}$ needs to hold all the variables and it suffices for $x_{i1}$ to only have a duplicate of $v_i$, i.e.
\bqn
x_{ir}:= \begin{cases}
(S_{i0}^{(x)}, \ell_{i0}^{(x)}, v_{i0}^{(x)}, s_{i0}^{(x)}) & r=0\\
(v_{i1}^{(x)}) & r=1 
\end{cases} .
\eqn
As a result, $x_{ir}$, $y_{ii}$ in \eqref{eq:distOPFu:eropf::con1} and $x_{i0}$, $y_{ij}$ in \eqref{eq:distOPFu:eropf::con2} do not consist of the same components. Here, we abuse notations in both \eqref{eq:distOPFu:eropf::con1} and \eqref{eq:distOPFu:eropf::con2}, which are composed of components that appear in both items, i.e.
\begin{align*}
&x_{i0}-y_{ij}\\
:=&
\begin{cases}
(S^{(x)}_{i0}-S^{(y)}_{ii},\ell^{(x)}_{i0}-\ell^{(y)}_{ii},v^{(x)}_{i0}-v^{(y)}_{ii},s^{(x)}_{i0}-s^{(y)}_{ii}) & \hspace{-0.1in} j=i\\
(S^{(x)}_{i0}-S_{iA_i}^{(y)}, \ell^{(x)}_{i0}-\ell_{iA_i}^{(y)}) & \hspace{-0.1in} j=A_i\\
(v^{(x)}_{i0}-v_{ij}^{(y)}) & \hspace{-0.1in} j\in C_i
\end{cases}\\
&x_{ir} - y_{ii}:=
\begin{cases}
(v^{(x)}_{i1}-v^{(y)}_{ii}) & r=1
\end{cases}.
\end{align*}

Let $\lambda$ denote the Lagrangian multiplier for \eqref{eq:distOPFu:eropf::con1} and $\mu$ the Lagrangian multiplier for \eqref{eq:distOPFu:eropf::con1}. The detailed mapping between constraints and those multipliers are illustrated in Table \ref{tab:distOPFu::muliplier}.

\begin{table}
\caption{Multipliers associated with constraints \eqref{eq:distOPFu:eropf::con1}-\eqref{eq:distOPFu:eropf::con2}}
\begin{center}
\begin{tabular}{|c|c||c|c|}
\hline
$\lambda_{i1}$: & $v_{i1}^{(x)}=v_{ii}^{(y)}$ &  & \\
\hline
$\mu^{(1)}_{ii}$: & $S_{i0}^{(x)}=S_{ii}^{(y)}$ & $\mu^{(2)}_{ii}$: & $\ell_{i0}^{(x)}=\ell_{ii}^{(y)}$\\
\hline
$\mu^{(3)}_{ii}$: & $v_{i0}^{(x)}=v_{ii}^{(y)}$ & $\mu^{(4)}_{ii}$: & $s_{i0}^{(x)}=s_{ii}^{(y)}$\\
\hline
$\mu^{(1)}_{iA_i}$: & $S_{iA_i}^{(x)}=S_{iA_i}^{(y)}$ & $\mu^{(2)}_{iA_i}$: & $\ell_{i}^{(x)}=\ell_{iA_i}^{(y)}$ \\
\hline
$\mu_{ij}$: & $v_{i}^{(x)}=v_{ij}^{(y)}$ & &\\ 
\hline
\end{tabular}
\end{center}
\label{tab:distOPFu::muliplier}
\end{table}

Next, we will derive the efficient solution for the subproblems in the $x$-update. For notational convenience, we will skip the iteration number $k$ on the variables. In the $x$-update, there are $2$ subproblems \eqref{eq:xupdatenode} associated with each bus $i$. The first problem, which updates $x_{i0}$, can be written explicitly as:
\begin{subequations}\label{eq:distOPFu::zagent_a}
\bq
\min && H_{i0}(x_{i0})  \label{eq:distOPFu::zagent_a1}\\
\mathrm{over} && x_{i0}=\{v_{i0}^{(x)}, \ell_{i0}^{(x)}, S_{i0}^{(x)}, s_{i0}^{(x)}\} \\
\mathrm{s.t.}&&    \begin{pmatrix}
v_{i0}^{(x)} & S_{i0}^{(x)}\\
(S_{i0}^{(x)})^H & \ell_{i0}^{(x)}
\end{pmatrix}\in \mathbb{S}_+  \label{eq:distOPFu::zagent_a2} \\
&& (s_{i0}^\phi)^{(x)}\in \mathI_i^\phi \qquad \phi\in\Phi_i, \label{eq:distOPFu::zagent_a3}
\eq
\end{subequations}
where $H_{i0}(x_{i0})$ is defined in \eqref{eq:admm:opt4:hix} and for our application, $\|x_{i0}-y_{ij}\|_{M_{ij}}^2$ is chosen to be
\begin{align}\label{eq:mijdef}
&\|x_{i0}-y_{ij}\|_{M_{ij}}^2=\\
&\begin{cases}
(2|C_i|+3)\|S_{i0}^{(x)}-S_{ii}^{(y)}\|_2^2+ \|s_{i0}^{(x)}-s_{ii}^{(y)}\|_2^2 &\\
\quad+2\|v_{i0}^{(x)}-v_{ii}^{(y)}\|_2^2+ (|C_i|+1)\|\ell_{i0}^{(x)}-\ell_{ii}^{(y)}\|_2^2& j=i\\
\|S_{i0}^{(x)}-S_{iA_i}^{(y)}\|_2^2+\|\ell_{i,A_i}^{(x)}-\ell_i^{(y)}\|_2^2 & j = A_i \\
\|x_{i0}-y_{ij}\|_{2}^2 & j\in C_i
\end{cases}\nonumber
\end{align}
By using \eqref{eq:mijdef}, $H_i^{(1)}(S_{i0}^{(x)}, \ell_{i0}^{(x)}, v_{i0}^{(x)})$, which is defined below, can be written as the Euclidean distance of two Hermitian matrix, which is one of the key reasons that lead to our efficient solution. Therefore, $H_{i0}(x_{i0})$ can be further decomposed as 
\begin{align}\label{eq:distOPFu::zupdate_square}
&H_{i0}(x_{i0})\\
=&f_i(x_{i0})+\sum_{j\in N_i}\left(\langle\mu_{ij},x_{i0}\rangle +\frac{\rho}{2}\|x_{i0}-y_{ij}\|_{M_{ij}}^2\right)\nonumber\\
=&\frac{\rho(|C_i|+2)}{2} H_i^{(1)}(S_{i0}^{(x)}, \ell_{i0}^{(x)}, v_{i0}^{(x)}) + H_i^{(2)}(s_{i0}^{(x)}) +\text{constant}, \nonumber
\end{align}
where 
\bqn
H_i^{(1)}(S_{i0}^{(x)}, \ell_{i0}^{(x)}, v_{i0}^{(x)})\hspace{-0.1in}&=&\hspace{-0.1in} \left\|
\begin{pmatrix}
v_{i0}^{(x)} & S_{i0}^{(x)} \\
(S_{i0}^{(x)})^H &  \ell_{i0}^{(x)}
\end{pmatrix}
-
\begin{pmatrix}
\hat v_i & \hat S_i\\
\hat S_i^H &  \hat \ell_i
\end{pmatrix}
\right\|_2^2\\
H_i^{(2)}(s_{i0}^{(x)})\hspace{-0.1in}&=&\hspace{-0.1in} f_i(s_{i0}^{(x)})+\frac{\rho}{2}\|s_{i0}^{(x)}-\hat s_i\|_2^2.
\eqn
The last step in \eqref{eq:distOPFu::zupdate_square} is obtained using square completion and the variables labeled with hat are some constants. 

Hence, the objective \eqref{eq:distOPFu::zagent_a1} in \eqref{eq:distOPFu::zagent_a} can be decomposed into two parts, where the first part $H^{(1)}(S_{i0}^{(x)}, \ell_{i0}^{(x)}, v_{i0}^{(x)})$ involves variables $(S_{i0}^{(x)}, \ell_{i0}^{(x)}, v_{i0}^{(x)})$ and the second part $H^{(2)}(s_{i0}^{(x)})$ involves $s_{i0}^{(x)}$. Note that the constraint \eqref{eq:distOPFu::zagent_a2}--\eqref{eq:distOPFu::zagent_a3} can also be separated into two independent constraints. Variables $(S_{i0}^{(x)}, \ell_{i0}^{(x)}, v_{i0}^{(x)})$ only depend on \eqref{eq:distOPFu::zagent_a2} and $s_{i0}^{(x)}$ depends on \eqref{eq:distOPFu::zagent_a3}. Then  \eqref{eq:distOPFu::zagent_a} can be decomposed into two subproblems, where the first one \eqref{eq:distOPFu::zagent1} solves the optimal $(S_{i0}^{(x)}, \ell_{i0}^{(x)}, v_{i0}^{(x)})$ and the second one \eqref{eq:distOPFu::zagent2} solves the optimal $s_{i0}^{(x)}$. The first subproblem can be written explicitly as 
\bq
\min &&H_i^{(1)}(S_{i0}^{(x)}, \ell_{i0}^{(x)}, v_{i0}^{(x)}) \nonumber\\
\mathrm{over} && S_{i0}^{(x)}, \ell_{i0}^{(x)}, v_{i0}^{(x)} \label{eq:distOPFu::zagent1}\\
\mathrm{s.t.}&&     \begin{pmatrix}
v_{i0}^{(x)} & S_{i0}^{(x)} \nonumber\\
(S_{i0}^{(x)})^H & \ell_{i0}^{(x)}
\end{pmatrix}\in \mathbb{S}_+ ,  \nonumber
\eq
which can be solved using eigen-decomposition of a $6\times 6$ matrix via the following theorem.
\begin{theorem}\label{thm:distOPFu::1}
Suppose $W\in\mathbb{S}^n$ and denote $X(W):=\arg\min_{X\in \mathbb{S}_+}\|X-W\|_2^2$. Then $X(W)=\sum_{i:\lambda_i>0}\lambda_iu_iu_i^H$, where $\lambda_i,u_i$ are the $i_{\text{th}}$ eigenvalue and orthonormal eigenvector of matrix $W$, respectively. 
\end{theorem}
\begin{IEEEproof}
The proof is in Appendix \ref{app:distOPFu::1}.
\end{IEEEproof}

Denote
\bqn
W:=
\begin{pmatrix}
\hat v_i & \hat S_i\\
\hat S_i^H &  \hat \ell_i
\end{pmatrix}
\text{ and }
X:=
\begin{pmatrix}
v_{i0}^{(x)} & S_{i0}^{(x)} \nonumber\\
(S_{i0}^{(x)})^H & \ell_{i0}^{(x)}
\end{pmatrix}.
\eqn
Then \eqref{eq:distOPFu::zagent1} can be written abbreviately as 
\bqn
\min_{X}\|X-W\|_2^2 \quad\text{s.t. }X \in \mathbb{S}_+,
\eqn
which can be solved efficiently using eigen-decomposition based on Theorem \ref{thm:distOPFu::1}. The second problem is 
\bq
\min && f_i(s_{i0}^{(x)})+\frac{\rho}{2}\|s_{i0}^{(x)}-\hat s_i\|_2^2 \nonumber\\
\mathrm{over} && s_{i0}^{(x)} \in\mathI_i^\phi  \qquad \phi\in\Phi_i . \label{eq:distOPFu::zagent2}
\eq
Recall that if $ f_i(s_{i0}^{(x)})=\sum_{\phi\in\Phi_i}f_i^{\phi}((s_{i0}^{\phi})^{(x)})$, then both the objective and constraint are separable for each phase $\phi\in\Phi_i$. Therefore, \eqref{eq:distOPFu::zagent2} can be further decomposed into $|\Phi_i|$ number of subproblems as below.
\bq
\min && f_i^{\phi}((s_{i0}^{\phi})^{(x)}) \nonumber\\
\mathrm{over} && (s_{i0}^\phi)^{(x)} \in\mathI_i^\phi , \label{eq:distOPFu::zagent21}
\eq
which takes the same form as of \eqref{eq:distOPFu:closedformsol}  in Theorem \ref{thm:distOPFu:closedformsol} and thus can be solved with closed form solution based on the assumptions. 

For the problem \eqref{eq:distOPFu::zagent_b} that updates $x_{i1}$, which consists of only one component $v_{i1}^{(x)}$, it can be written explicitly as
\bq
\min && H_{i1}(x_{i1})   \nonumber\\
\mathrm{over} && x_{i1}=\{v_{i1}^{(x)}\} \label{eq:distOPFu::zagent_b}\\
\mathrm{s.t.}&&   \underline v_i^\phi \leq (v_{i1}^{\phi\phi})^{(x)} \leq \overline v_i^\phi  \qquad \phi\in\Phi_i ,\nonumber 
\eq
where $H_{i1}(x_{i1})$ is defined in \eqref{eq:admm:opt4:hix} and for our application, $\|x_{ir}-y_{ii}\|_{\Lambda_{ir}}^2$ is chosen to be
\bqn
\|x_{ir}-y_{ii}\|_{\Lambda_{ir}}^2=\|x_{ir}-y_{ii}\|_{2}^2
\eqn
Then the closed form solution is given as:
\bqn
(v_{i1}^{\phi_1\phi_2})^{(x)}=\begin{cases}
\left[\frac{\lambda_{i1}^{\phi_1\phi_2}}{\rho}+(v_{ii}^{\phi_1\phi_2})^{(y)}\right]_{\underline v_i^{\phi_1}}^{\overline v_i^{\phi_1}} & \phi_1=\phi_2\\
\frac{\lambda_{i1}^{\phi_1\phi_2}}{\rho}+(v_{ii}^{\phi_1\phi_2})^{(y)} & \phi_1\neq \phi_2
\end{cases}.
\eqn
To summarize, the subproblems in the $x$-update for each bus $i$ can be solved either through a closed form solution or a eigen-decomposition of a $6\times6$ matrix, which proves Theorem \ref{thm:distOPFu:closedformsol}.

In the $y$-update, the subproblem solved by each node $i$ takes the form of \eqref{eq:admm::opt4::zupdate1} and can be written explicitly as 
\begin{align}
\min \ &G_i(\{y_{ji}\mid j\in N_i\})\nonumber\\
\mathrm{over} \ & \{y_{ji}\mid j\in N_i\} \label{eq:distOPFu::xagent}\\
\mathrm{s.t.} \ & \mathP_i(v^{(y)}_{A_ii})=v^{(y)}_{ii}-z_i(S^{(y)}_{ii})^H-S^{(y)}_{ii}z_i^H+z_i\ell^{(y)}_{ii}z_i^H\nonumber  \\
\ &s^{(y)}_{ii}=-\text{diag}\left(\sum_{i\in C_i}\mathP_i(S^{(y)}_{ji}-z_j\ell^{(y)}_{ji})-S^{(y)}_{ii}\right),  \nonumber
\end{align}
which has a closed form solution given in \eqref{eq:admm::opt2::zsol} and we do not reiterate here. 

Finally, we specify the initialization and stopping criteria for the algorithm. Similar to the algorithm for balanced networks, a good initialization usually reduces the number of iterations for convergence. We use the following initialization suggested by our empirical results. We first initialize the auxiliary variables $\{V_i\mid i\in\mathN\}$ and $\{I_i\mid i\in\mathE\}$, which represent the complex nodal voltage and branch current, respectively. Then we use these auxiliary variables to initialize the variables in \eqref{eq:distOPFu::eropf}. Intuitively, the above initialization procedure can be interpreted as finding a solution assuming zero impedance on all the lines. The procedure is formally stated in Algorithm \ref{alg:distOPFu::initialize}.

  \begin{algorithm}
 \caption{Initialization of the Algorithm}
 \label{alg:distOPFu::initialize}
  \begin{algorithmic}[1]
 \State $V_i^{a}=1$, $V_i^{b}=e^{-\frac{2}{3}\pi}$, $V_i^{c}=e^{\frac{2}{3}\pi}$ for $i\in\mathN$
 \State Initialize $s_i^\phi$ using any point in the injection region $\mathI_i^\phi$ for $i\in\mathN$
 \State Initialize $\{I_i^\phi\mid \phi\in\Phi_i \ i\in\mathN\}$ by calling DFS($0$,$\phi$) for $\phi\in\Phi_i$
 \State $v_{i0}^{(x)}=V_iV_i^H$, $\ell_{i0}^{(x)}=I_iI_i^H$, $S_{i0}^{(x)}=V_iI_i^H$ and $s_{i0}^{(x)}=s_i$ for $i\in\mathN$
 \State $y_{ij}=x_{i0}$ for $j\in N_i$ and $i\in\mathN$ 
 \State $x_{i1}=x_{i0}$ for $i\in \mathN$
 \Statex
\Function{DFS}{$i$,$\phi$}
    \State $I_i^\phi=(\frac{s_i^{\phi}}{V_i^\phi})^*$ 
    \For{$j\in C_i$}
        \State $I_i^{\phi}+=\mathrm{DFS}(j,\phi)$ 
      \EndFor
    \State \Return $I_i^\phi$
\EndFunction
 \end{algorithmic}
 \end{algorithm}

For the stopping criteria, there is no general rule for ADMM based algorithm and it usually hinges on the particular problem. In \cite{boyd2011distributed}, it is argued that a reasonable stopping criteria is that both the primal residual $r^k$ defined in \eqref{eq:pfeasible} and the dual residual $s^k$ defined in \eqref{eq:dfeasible} are below $10^{-4}\sqrt{|\mathN|}$. We adopt this criteria and the empirical results show that the solution is accurate enough. The pseudo code for the complete algorithm is summarized in Algorithm \ref{alg:distOPFu::alg}.

  \begin{algorithm}
 \caption{Distributed OPF algorithm on Unbalanced Radial Networks}
 \label{alg:distOPFu::alg}
  \begin{algorithmic}[1]
 \State \algorithmicrequire  network $\mathG(\mathN,\mathE)$, power injection region $\mathI_i$, voltage region $(\underline v_i,\overline v_i)$, line impedance $z_i$ for $i\in\mathN$.
 \State \algorithmicensure  voltage $v$, power injection $s$
 \Statex
 \State Initialize the $x$ and $y$ variables using Algorithm \ref{alg:distOPFu::initialize}.
 \While{$r^k>10^{-4}\sqrt{|\mathN|}$ \textbf{or}  $s^k>10^{-4}\sqrt{|\mathN|}$ }
 \State  In the $x$-update, each agent $i$ solves both \eqref{eq:distOPFu::zagent_a} and \eqref{eq:distOPFu::zagent_b} to update $x_{i0}$ and $x_{i1}$.
 \State In the $y$-update, each agent $i$ solves \eqref{eq:distOPFu::xagent} to update $y_{ji}$ for $j\in N_i$.
 \State In the multiplier update, update $\lambda,\mu$ by  \eqref{eq:mupdate}.
 \EndWhile
 \end{algorithmic}
 \end{algorithm}


\section{Case Study}\label{sec:case}

In this section, we first demonstrate the scalability of the distributed algorithm proposed in section \ref{sec:dalg} by testing it on the standard IEEE test feeders \cite{kersting1991radial}. To show the efficiency of the proposed algorithm, we also compare the computation time of solving the subproblems in both the $x$ and $y$-update with off-the-shelf solver (CVX). Second, we run the proposed algorithm on networks of different topology to understand the factors that affect the convergence rate. The algorithm is implemented in Python and run on a Macbook pro 2014 with i5 dual core processor. 

\subsection{Simulations on IEEE test feeders}\label{sec:distOPFu::case1}

We test the proposed algorithm on the IEEE 13, 34, 37, 123 bus distribution systems. All the networks have unbalanced three phase. The substation is modeled as a fixed voltage bus ($1$ p.u.) with infinite power injection capability. The other buses are modeled as load buses whose voltage magnitude at each phase can vary within $[0.95,1.05]$ p.u. and power injections are specified in the test feeder. There is no controllable device in the original IEEE test feeders, and hence the OPF problem degenerates to a power flow problem, which is easy solve. To demonstrate the effectiveness of the algorithm, we replace all the capacitors with inverters, whose reactive power injection ranges from $0$ to the maximum ratings specified by the original capacitors. The objective is to minimize power loss across the network, namely $f_i^{\phi}(s_i^\phi)=p_i^\phi$ for $\phi\in\Phi_i$ and $i\in\mathN$.

We mainly focus on the time of convergence (ToC) for the proposed distributed algorithm. The algorithm is run on a single machine. To roughly estimate the ToC (excluding communication overhead) if the algorithm is run on multiple machines in a distributed manner, we divide the total time by the number of buses. 

\begin{table}
\caption{Statistics of different networks}
\begin{center}
\begin{tabular}{|c|c|c|c|c|c|}
\hline
Network  &Diameter& Iteration & Total Time(s) & Avg time(s)\\
\hline
IEEE 13Bus &6& 289 & 17.11 & 1.32\\
\hline
IEEE 34Bus &20&547 & 78.34 & 2.30\\ 
\hline
IEEE 37Bus &16&440 &75.67 & 2.05\\ 
\hline
IEEE 123Bus &30&608 &306.3 & 2.49\\ 
\hline
\end{tabular}
\end{center}
\label{tab:statistics}
\end{table}

In Table \ref{tab:statistics}, we record the number of iterations to converge, total computation time required to run on a single machine and the average time required for each node if the algorithm is run on multiple machines excluding communication overhead. From the simulation results, the proposed algorithm converges within $2.5$ second for all the standard IEEE test networks if the algorithm is run in a distributed manner. 

Moreover, we show the advantage of using the proposed algorithm by comparing the computation time of solving the subproblems between off-the-shelf solver (CVX \cite{grant2008cvx}) and our algorithm. In particular, we compare the average computation time of solving the subproblem in both the $x$ and $y$ update. In the $x$-update, the average time required to solve the subproblem \eqref{eq:distOPFu::xagent} is $9.8\times10^{-5}$s for our algorithm but $0.13$s for CVX.  In the $y$-update, the average time required to solve the subproblems \eqref{eq:distOPFu::zagent_a}--\eqref{eq:distOPFu::zagent_b} are $3.7\times10^{-3}$s for our algorithm but $0.45$s for CVX. Thus, each ADMM iteration takes about $3.8\times 10^{-3}$s for our algorithm but $5.8\times 10^{-1}$s for using iterative algorithm, a more than 100x speedup.

\subsection{Impact of Network Topology}

In section \ref{sec:distOPFu::case1}, we demonstrate that the proposed distributed algorithm can dramatically reduce the computation time within each iteration. The time of convergence (ToC) is determined by both the computation time required within each iteration and the number of iterations. In this subsection, we study the number of iterations, namely rate of convergence. 

Rate of convergence is determined by many different factors. Here, we only consider the rate of convergence from two factors, network size $N$, and diameter $D$, i.e. given the termination criteria in Algorithm \ref{alg:distOPFu::alg}, the impact of network size and diameter on the number of iterations. The impact from other factors, e.g. form of objective function and constraints, is beyond the scope of this paper. 

\begin{figure}
\centering
\subfloat[Line network] {
	\includegraphics[scale=0.35]{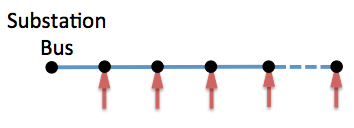}
	\label{fig:line}
	}
\subfloat[Fat tree network]{
	\includegraphics[scale=0.35]{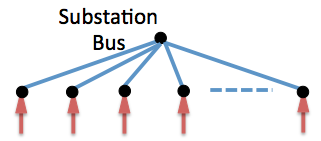}
	\label{fig:fattree}
	}
\caption{Topologies for line and fat tree networks.}
\label{fig:network}
\end{figure}

To illustrate the impact of network size $N$ and diameter $D$ on the rate of convergence, we simulate the algorithm on two extreme cases: 1) Line network in Fig. \ref{fig:line}, whose diameter is the largest given the network size, and 2) Fat tree network in Fig. \ref{fig:fattree}, whose diameter is the smallest given the network size. In Table \ref{tab:statistics2}, we record the number of iterations for both line and fat tree network of different sizes. For the line network, the number of iterations increases notably as the size increases. For the fat tree network, the trend is less obvious compared to line network. It means that the network diameter has a stronger impact than the network size on the rate of convergence.

\begin{table}
\caption{Statistics of line and fat tree networks}
\begin{center}
\begin{tabular}{|c|c|c|}
\hline
Size & $\#$ of iterations (Line) & $\#$ of iterations (Fat tree) \\
\hline
$5$ & $57$ & $61$ \\
\hline
$10$ & $253$  & $111$\\
\hline
$15$ & $414$ & $156$\\
\hline
$20$ & $579$ & $197$\\
\hline 
$25$ & $646$ & $238$\\
\hline 
$30$ & $821$ & $272$\\
\hline
$35$ & $1353$ & $304$\\ 
\hline
$40$ & $2032$ & $337$\\
\hline
$45$ & $2026$ & $358$\\
\hline
$50$ & $6061$ & $389$\\
\hline
\end{tabular}
\end{center}
\label{tab:statistics2}
\end{table}

\section{Conclusion}\label{sec:conclusion}
In this paper, we have developed a distributed algorithm for optimal power flow problem on unbalanced distribution system based on alternating direction method of multiplier. We have derived an efficient solution for the subproblem solved by each agent thus significantly reducing the computation time. Preliminary simulation shows that the algorithm is scalable to all IEEE test distribution systems.

\bibliographystyle{IEEEtran}
\bibliography{texcode/reference}{}

\appendices
\section{Proof of Theorem \ref{thm:distOPFu::1}}\label{app:distOPFu::1}
Let $\Lambda_W:=$diag$(\lambda_i,1\leq i\leq n)$ denote the diagonal matrix consisting of the eigenvalues of matrix $W$. 
Let $U:=(u_i,1\leq i\leq n)$ denote the unitary matrix. Since $W\in\mathbb{H}^n$, $U^{-1}=U^H$ and $W=U\Lambda_WU^H$.
Then
\bqn
\|X-W\|_2^2&=&tr((X-W)^H(X-W))\\
&=&tr((X-W)(X-W))\\
&=&tr(U^H(X-W)UU^H(X-W)U)\\
&=&tr((U^HXU-\Lambda_W)(U^HXU-\Lambda_W)).
\eqn
Denote $\hat X:=U^HXU=(\hat x_{i,j},i,j\in[1,n])$, note that $\hat X\in\mathbb{S}_+$ since $X\in\mathbb{S}_+$. Then
\bq
\|X-W\|_2^2&=& \sum_{i=1}^n(\hat x_{ii}-\lambda_i)^2+\sum_{i\neq j}|\hat x_{ij}|^2\\
&\geq& \sum_{i=1}^n(\hat x_{ii}-\lambda_i)^2\\
&\geq& \sum_{i:\lambda_i\leq 0}\lambda_i^2, \label{eq:fbound}
\eq
where the last inequality follows from $\hat x_{ii}>0$ because $\hat X\in\mathbb{S}_+$. The equality in \eqref{eq:fbound} can be obtained by letting 
\bqn
\hat x_{ij}:=\left\{
                \begin{array}{ll}
                  \lambda_i & i=j, \ \ \lambda_i>0, \\
                  0 & \text{ otherwise}
                \end{array},
              \right.
\eqn
which means $X(W)=U\hat X U^H=\sum_{i:\lambda_i>0}\lambda_iu_iu_i^H$.

\section{Solution Procedure for Problem \eqref{eq:distOPFu:closedformsol}.}\label{app:distOPFb::solver2}

We assume $f_i^\phi\left(s\right):=\frac{\alpha_i}{2} p^2+\beta_i p$ $(\alpha_i, \beta_i\geq 0)$ and derive a closed form solution to \eqref{eq:distOPFu:closedformsol}

\subsection{$\mathI_i$ takes the form of \eqref{eq:opf::unbalance::I2}}
In this case, \eqref{eq:distOPFu:closedformsol} takes the following form:
\bqn
\min_{p,q} && \frac{a_1}{2}p^2+b_1p+\frac{a_2}{2}q^2+b_2q\\
\text{s.t. } && \underline p_i\leq p\leq \overline p_i \\
&& \underline q_i\leq q\leq \overline q_i ,
\eqn
where $a_1,a_2>0$ and $b_1,b_2$ are constants. Then the closed form solution is 
\bqn
p=\left[-\frac{b_1}{a_1}\right]_{\underline p_i}^{\overline p_i} \quad
q=\left[\hat -\frac{2_1}{a_2}\right]_{\underline q_i}^{\overline q_i} ,
\eqn
where $[x]_a^b:=\min\{a,\max\{x,b\}\}$.

\subsection{$\mathI_i$ takes the form of \eqref{eq:opf::unbalance::I1}}
The optimization problem \eqref{eq:distOPFu:closedformsol} takes the following form:
\begin{subequations}\label{eq:app2}
\bq
\min_{p,q} && \frac{a_1}{2}p^2+b_1p+\frac{a_2}{2}q^2+b_2q\\
\text{s.t. } && p^2+q^2\leq c^2 \label{eq:app2:1}\\
&& p \geq 0 , \label{eq:app2:2}
\eq
\end{subequations}
where $a_1,a_2,c>0$ $,b_1,b_2$ are constants. The solutions to \eqref{eq:app2} are given as below.
{\flushleft\bf Case 1}: $b_1\geq 0$:
\bqn
p^*=0 \qquad q^*=\left[-\frac{b_2}{a_2}\right]_{-c}^{c}.
\eqn
{\flushleft\bf Case 2}: $b_1< 0$ and $\frac{b^2_1}{a^2_1}+\frac{b^2_2}{a_2^2}\leq c^2$:
\bqn
p^*=-\frac{b_1}{a_1} \qquad q^*=-\frac{b_2}{a_2}.
\eqn
{\flushleft\bf Case 3}: $b_1< 0$ and $\frac{b^2_1}{a^2_1}+\frac{b^2_2}{a_2^2}> c^2$:\\
First solve the following equation in terms of variable $\lambda$:
\bq
b_1^2(a_2+2\lambda)^2+b_2^2(a_1+2\lambda)^2=(a_1+2\lambda)^2(a_2+2\lambda)^2, \label{eq:app3:1}
\eq
which is a polynomial with degree of $4$ and has closed form expression. There are four solutions to \eqref{eq:app3:1}, but there is only one strictly positive $\lambda^*$, which can be proved via the KKT conditions of \eqref{eq:app2}. Then we can recover $p^*,q^*$ from $\lambda^*$ using the following equations:
\bqn
 p^*=-\frac{b_1}{a_1+2\lambda^*} \quad \text{ and }\quad q^*=-\frac{b_2}{a_2+2\lambda^*}.
\eqn

The above procedure to solve \eqref{eq:app2} is derived from standard applications of the KKT conditions of \eqref{eq:app2}. For brevity, we skip the proof here.



\end{document}